\documentclass[12pt]{article}

\usepackage{amssymb}
\makeatletter

\@addtoreset{equation}{section}
\def\section{\@startsection {section}{1}{\z@}{-3.5ex plus -1ex minus
 -.2ex}{2.3ex plus .2ex}{\large\bf}}
\def\subsection{\@startsection{subsection}{2}{\z@}{-3.25ex plus%
 -1ex minus -.2ex}{1.5ex plus .2ex}{\sc}}
 
\advance \voffset by -0.8in
\advance \hoffset by -0.6in
\def\ad{\mbox{Ad}}
\def\cd{\!\cdot\!}
\def\calC{{\mathcal C}}
\def\ba{{\mbox{\boldmath $a$}}}
\def\bj{{\mbox{\boldmath $j$}}}
\def\bk{{\mbox{\boldmath $k$}}}

\def\bd{{\mbox{\boldmath $d$}}}

\def\bn{{\mbox{\boldmath $n$}}}
\def\bob{{\mbox{\boldmath $b$}}}
\newcommand{\gothg}{\mathfrak g }
\newcommand{\gothh}{\mathfrak h }
\newcommand{\NN}{\mathbb{N}}
\newcommand{\ZZ}{\mathbb{Z}}

\newcommand{\RR}{\mathbb{R}}
\newcommand{\CC}{\mathbb{C}}

\def\bea{\begin{eqnarray}}
\def\eea{\end{eqnarray}}

\begin{document}

\begin{flushright}
MS-00-007\\
math.QA/0006228
\end{flushright}
\vspace{0.5cm}
\begin{center}
\baselineskip 24 pt
{\LARGE Combinatorial quantisation of  Euclidean gravity  } \\
 {\LARGE in three  dimensions}

\baselineskip 18pt
\parskip 7pt

\vspace{1cm}

{\large B.J. Schroers\footnote{address after 1.09.2000:
Department of Mathematics,
Heriot-Watt University,
Edinburgh EH14 4AS, United Kingdom}  \\
Department of  Mathematics and Statistics, University of Edinburgh\\
King's Buildings, Mayfield Road, Edinburgh EH9 3JZ\\
United Kingdom }\\
 {\tt bernd@maths.ed.ac.uk}

\vspace{0.5cm}

{\large 10 June 2000}

\end{center}

\vspace{0.5cm}

\begin{abstract}
In the Chern-Simons formulation of Einstein gravity in 2+1 dimensions
the phase space of gravity is the moduli space of flat $G$-connections,
where $G$ is a typically  non-compact Lie group which depends on the signature 
of space-time and the cosmological constant. For Euclidean signature 
and vanishing cosmological constant, $G$ is the three-dimensional Euclidean
group. For this case the Poisson structure of the moduli space
is given explicitly in terms of a classical $r$-matrix.
It is shown that the quantum $R$-matrix of the quantum double 
 $D(SU(2))$ provides a quantisation of that Poisson structure.
\end{abstract}

\vspace{0.5cm}
\centerline{MSC 17B37, 81R50, 81S10, 83C45}

\section{Introduction}

The primary goal of this paper is to indicate how some of the  quantisation
techniques developed in the quantisation of Chern-Simons theory
with a compact gauge group can be extended and applied to
the quantisation of  three-dimensional gravity.
One expects this to be possible because gravity in three dimensions
can be re-formulated as a Chern-Simons theory. 
The gauge group of the Chern-Simons theory,
however, depends on  the  cosmological constant and the signature of 
space-time and is non-compact in almost all cases. 
 Here we shall follow the combinatorial 
quantisation program developed for Chern-Simons theory
with a compact gauge group  by Alekseev,
Grosse  and Schomerus, see \cite{AGSI}\cite{AGSII}
and also \cite{AS}. We shall show how  to implement the main steps of that 
program in one particular case, namely Euclidean gravity without
cosmological constant. While our results provide key stepping stones
on a   promising path to a 
 quantisation of three dimensional gravity, a number of issues - 
both physical and mathematical - 
are only raised but not settled here.
A secondary purpose of this paper is  to 
advertise these issues  to   
 mathematicians and physicists
with an active interest in the geometry and quantisation of 
 moduli spaces  of flat connections.

The paper naturally falls into two halves. The first half, consisting of
sects. 2 and 3, is a review  of the Chern-Simons formulation
 of three-dimensional gravity. The precise relation between this
formulation and the original Einstein formulation of gravity
is a bone of contention in the literature. We indicate some 
of the issues but do not enter deeply into their discussion. 
Our view is that the Chern-Simons formulation offers a promising
avenue towards quantising three-dimensional gravity, 
 which is worth pursuing because
of the significance of the goal and the possibility of concrete
results. At the end of sect. 3 the problem of quantising 
three-dimensional gravity  will have been translated into the 
problem of quantising certain moduli spaces.

The new results of this paper are contained in sects. 4 and 5,
which deal with the quantisation problem  for 
Euclidean gravity with vanishing cosmological constant.
In that 
case the gauge group of the Chern-Simons formulation is the (double
cover of the) Euclidean group $ISO(3)$ in three dimensions. If space-time
is a direct product of time and a two-dimensional space $\Sigma$, 
the phase space of gravity can  be identified with the space of 
flat $ISO(3)$ bundles on $\Sigma$. 
In the framework of combinatorial  quantisation, 
 the starting point of the quantisation is the 
Fock-Rosly  description of  the Poisson  structure of the phase space in terms
of a classical $r$-matrix, i.e. a solution of the classical
Yang-Baxter equation \cite{FR}. 
The key step in the quantisation of the phase space 
is the identification of a solution of the quantum Yang-Baxter
equations which reduces to the classical $r$-matrix in a suitable
limit.  Here we shall give the  $r$-matrix  relevant for $ISO(3)$
Chern-Simons theory and show that it is the limit of the universal
$R$-matrix of the quantum double $D(SU(2))$. The pivotal role of the 
quantum double $D(SU(2))$ in the quantisation of $ISO(3)$ Chern-Simons
theory was discoverd by Bais and Muller in
\cite{BM}. In our sect. 5 we shall show that $D(SU(2))$ is
a deformation  of the group algebra of $ISO(3)$  and 
thereby resolve a  question posed in \cite{BM}. 

In this paper we  restrict attention to 
Euclidean gravity and vanishing cosmological constant.
 The physically more interesting case 
of Lorentzian signature with vanishing cosmological constant 
can be dealt with in an analogous manner. However, 
some additional technical problems and numerous physical
implications call for a more detailed  discussion, which we give
in a separate paper \cite{BMS}. 
The inclusion of a non-vanishing cosmological constant
in the combinatorial quantisation of three-dimensional gravity poses 
a very interesting problem.
 This problem is addressed in   \cite{BR},  and we will
briefly comment on it from our viewpoint  at the end of this paper.

\section{The Chern-Simons formulation of gravity in three dimensions}

The possibility of writing general relativity
 in three  dimensions as a Chern-Simons
theory was first noticed in \cite{AT}. This observation opened up a new 
approach to gravity and in particular to its quantisation, which was
first systematically explored in \cite{Witten1}. Since then, a vast 
body of literature has been devoted to the subject. This section  
is an attempt to give the  briefest possible   summary of the Chern-Simons
formulation of three-dimensional Euclidean gravity. For more background
and references on three-dimensional gravity we refer the reader to
the recent book by Carlip \cite{Carlipbook} or the  review article
\cite{Carlipart}.

In  three dimensional gravity, space-time is a three-dimensional
manifold $M$.  In the following we shall only consider 
 space times of the form $M=\RR\times \Sigma$, where $\Sigma$
is an orientable two-dimensional manifold (``space''). 
A three-manifold  of that form 
is  orientable  and hence, by a classic theorem of  Stiefel, 
  parallelisable. Thus its tangent bundle is topologically
trivial.

 In Einstein's original formulation of general 
relativity,  the dynamical variable is a  metric $g$ on
$M$. For our purposes it
is essential  to adopt  Cartan's point of view, where the 
theory is formulated as a gauge theory. In this approach one
introduces an auxiliary 3-dimensional vector bundle $V$ with
 an inner product $(\,\, ,\,\,)$  and connection $\omega$,
metric with respect to $(\,\, ,\,\,)$.
The topological type of $V$ is that of the tangent bundle $TM$
of  $M$ (i.e. trivial in our case) and 
the structure group of $V$ is  $SO(3)$ (in the Lorentzian case it would be 
the Lorentz group $SO(2,1)$). 
Then there exists a   bundle map $TM \rightarrow V$
covering the identity. Such a bundle map provides an identification of 
$T_xM$ with the fibre $V_x$ of $V$ over $x\in M$, and can be thought
of as a V-valued one-form (soldering form or dreibein) $e$ on $M$. 
Choosing a basis $\{E_a\}$, $a=1,2,3$, 
of $V_x$, orthonormal with respect to $(\,\, ,\,\,)$
 and local coordinates $x_\mu$, $\mu =1,2,3$, around $x$,
we require   that the $3\times 3$  matrix $((e^a_\mu))$ 
defined by $e(\partial_\mu)
=\sum_{a=1}^3 e^a_\mu E_a$ is invertible.

To continue, we introduce generators $J_a$ of the Lie algebra
$so(3)$. They are normalised to satisfy
\bea
\label{spinrels}
[J_a,J_b]=\epsilon_{abc} J_c,
\eea
where $\epsilon_{abc}$ is the totally antisymmetric tensor in 
three dimensions, normalised so that $\epsilon_{123}=1$.  Here 
and in the following, repeated indices are summed on;
 since we are in the Euclidean situation, 
the position of indices  (upstairs or downstairs) is purely
for notational convenience.
The connection one-form $\omega$ can be expanded as $\omega=
 \omega_a J^a$. Similarly  the curvature two-form $F_\omega = d\omega +
\frac{1}{2}[\omega,\omega]$ can be expanded as $F_\omega=
F^a_\omega J_a$, with 
\bea
\label{spincurv}
F_\omega^a = d\omega^a  + \frac{1}{2} \epsilon^a_{\,bc} \omega^b\wedge 
\omega^c.
\eea
The Einstein-Hilbert action in  three dimension can be written as 
\bea 
\label{EHaction}
S_{EH}[\omega, e]= \int_M \,\, e_a \wedge F_\omega^a.
\eea
In Cartan's formulation, both the connection $\omega$ and the 
dreibein $e$ should be thought of as dynamical variables and 
varied independently. Variation with respect to $\omega$
yields the requirement that the connection  $\omega$ has 
vanishing torsion:
\bea
\label{notorsion}
D_\omega e_a = de_a + \frac{1}{2} \epsilon_{abc} \omega^b e^c =0.
\eea
This condition, which is imposed {\it a priori} in Einstein's
formulation of general relativity, is thus seen to be part of 
the equation of motion in Cartan's formulation.
 Variation with respect to $e$
yields the vanishing of the curvature tensor:
\bea
\label{Einstein}
F_{\omega}=0.
\eea
In   three  dimensions this is equivalent to the vanishing of the 
Ricci tensor, and thus to Einstein equations in  the 
absence of matter.

The geometrical background to Cartan's formulation is explained
beautifully  in the book \cite{Sharpe}. To explain the key ideas
in the present context we  choose a global frame $E_a$ and 
consider globally defined soldering forms $e_a$ on $M$. 
A Cartan connection may then  be defined as a  one-form
with values in the Lie algebra $iso(3)$ of the Euclidean group 
$ISO(3)=\RR^3\rtimes SO(3)$. 
Thus, if  we introduce translation generators $P_a$, $a=1,2,3$,
which satisfy 
\bea 
\label{transrels}
[J_a,P_b]=\epsilon_{abc} P^c, \quad
[P_a,P_b] = 0, 
\eea 
the Cartan connection can be written as 
\bea
\label{Cartan}
A = \omega_a J^a + e_a P^a.
\eea
The Cartan connection should be contrasted  with the usual Ehresmann
notion of a connection on a principal fibre bundle. 
While the Ehresmann connection is a one-form with values in the 
Lie-algebra of the structure group, the Cartan connection takes 
values in a bigger Lie algebra.
In the present case, it is a connection on a principal $SO(3)$
bundle, but takes values in  $iso(3)$. 
Cartan connections also have to satisfy a non-degeneracy condition,
which in the present context requires that the  soldering forms $e_a$
are nowhere vanishing.  Finally we note that the curvature of the 
Cartan connection 
\bea
\label{decomp}
F = (D_\omega e^a) P_a + (F_\omega^a) J_a
\eea
combines the curvature  and the torsion of the spin connection.

The Cartan framework allows one to translate  Riemannian (or Lorentzian)
geometry  into an equivalent gauge theory. However, the resulting
gauge theory is not of the standard type, and conditions have to be 
imposed on the connection. The crucial - and contentious - 
step in rewriting writing three dimensional 
gravity as a Chern-Simons theory is to drop these conditions and 
to interpret (\ref{Cartan}) as an Ehresmann 
connection of a bundle whose structure group is 
$ISO(3)$ . Thus, in particular the condition
of the invertibility of $e^a_\mu$ is dropped.  We will not enter into
the discussion of the merits and drawbacks of this approach. 
The advantages 
for  the quantisation are explained in Witten's original
paper \cite{Witten1}. For a recent, carefully argued 
 criticism, see \cite{Matschull}.

A final technical ingredient we need in order to establish the Chern-Simons
formulation is  special to  three dimensional space-times.
This is a non-degenerate, invariant
  bilinear form on the Lie algebra $iso(3)=
 \RR^3\rtimes so(3)$:
\bea
\label{inprod}
\langle J_a, P_a\rangle = \delta_{ab}, \quad  \langle J_a, J_b\rangle
= \langle P_a,P_b\rangle = 0.
\eea
Note that with respect to this inner product 
 both $so(3)$  and $\RR^3$ are maximally and 
totally null. 

Finally, we can write down  the Chern-Simons action on for the
connection
$A$ on  $M=\Sigma\times \RR$:
\bea
\label{CSaction}
S_{CS}[A] =\frac{1}{2} \int_M \langle A\wedge  dA\rangle 
 +\frac{2}{3}\langle A \wedge A \wedge A\rangle.
\eea 
A short calculation shows that this is equal to the Einstein-Hilbert
 action (\ref{EHaction}). Moreover, the equation of motion found 
by varying the action with respect to $A$ is 
\bea
\label{flat}
F=0.
\eea
Using the decomposition (\ref{decomp}) we thus reproduce
the condition of vanishing torsion and the  three dimensional
Einstein equations, as required.

So far we have only studied the Einstein equations in vacuum.
The introduction of matter in the form of point particles 
is physically desirable. Happily, it can be implemented 
in a mathematically elegant fashion in the Chern-Simons 
formulation. We refer the reader 
to \cite{BM} for a detailed discussion and further references,
 and only summarise
the salient points here. Particles are introduced by marking points
on the surface $\Sigma$ and coupling  the particle's phase space
to the phase space of the theory. The phase space of a 
 particle with Euclidean  mass $\mu$ and  spin $s$ is 
 a co-adjoint orbit
  ${\mathcal O}_{\mu s}$ 
of $ISO(3)$. To describe these orbits we  write $P_a^*$ and 
$J_a^*$ for the basis elements of $iso(3)^*$ dual to $P_a$ and $J_a$,
and we write an element  $\xi^*\in iso(3)^*$ as 
\bea
\label{dualphase}
\xi^*=p^a P_a^* + j^a J_a^*.
\eea
Using the inner product 
 (\ref{inprod}) we can identify $\xi^*$ with the element
\bea
\label{phaseco}
\xi=p^a J_a + j^a P_a 
\eea
in $iso(3)$.
Then $p^a$ should be thought of as the energy-momentum vector of the 
particle and $j^a$ as its generalised angular momentum.
The orbit ${\mathcal O}_{\mu s}$ consists of all $\xi^\star\in iso(3)^*$
satisfying  the mass-shell condition
$p_a p^a=\mu^2$ and the spin condition $p_aj^a=\mu s$.
As explained in detail in the book \cite{MR},
 ${\mathcal O}_{0 0}$ is a point, ${\mathcal O}_{0 s}$, with $s\neq 0$,  is 
a two-sphere of radius $s$ and ${\mathcal O}_{\mu s}$ for $\mu\neq 0$ 
 is diffeomorphic to $TS^2$. Co-adjoin torbits have a canonical
symplectic structure, often called the Kostant-Kirillov symplectic structure.
For the generic case $\mu\neq 0$ the corresponding Poisson brackets 
of the coordinate functiosn $j_a$ and $p_a$ are
\bea 
\label{kiri}
\{j_a,j_b\}=\epsilon_{abc}j_c, \qquad \{j_a,p_b\}=\epsilon_{abc}p_c.
\eea

In order to introduce $m$ particles with masses and 
spins $(\mu_1,s_1),...(\mu_m,s_m)$ we thus mark $m$ points $z_1, ...,z_m$
on $\Sigma$ and associate to each point $z_i$
a co-adjoint orbit ${\mathcal O}_{\mu_i s_i}$,$i=1,...m$. The coupling 
of the particle degrees of freedom to the gauge field via minimal
coupling is described in \cite{BM}. The upshot is that 
 we specify the kinematic state of each particle by  picking
elements $\xi^*_{(i)}=p_{(i)}^a P_a^* + j_{(i)}^a J_a^* \in 
{\mathcal O}_{\mu_i s_i}$.
The dual $iso(3)$ elements
$\xi_{(i)}=p_{(i)}^a J_a + j_{(i)}^a P_a $ then  act as the sources of 
curvature at  each of the marked points:
\bea
\label{gravwith}
F=\sum_{i=1}^m (p_{(i)}^a J_a+ j_{(i)}^a P_a)\delta(z-z_i).
\eea
Expanding the curvature term as in (\ref{decomp}) we find that 
the energy-momentum vectors of the particles act
 as  sources for  curvature and their
generalised angular momenta 
act as  sources of torsion, in agreement  with physical 
expectations.

The above discussion can be generalised to include a  non-vanishing 
cosmological constant $\lambda\in \RR$  and Lorentzian gravity. We refer the 
reader to \cite{Witten1} and \cite{Witten2} for details. The idea
is again to  
 combine the spin connection and the dreibein into a  Cartan
connection. The form of the Cartan connection
remains (\ref{Cartan})  but the 
Lie algebra structure of the space spanned by the generators 
$J_a$ and $P_a$ is modified as follows
\bea
\label{comrels}
[J_a,J_b]=\epsilon_{abc} J^c, \quad
[J_a,P_b]=\epsilon_{abc} P^c, \quad
[P_a,P_b]=\lambda \epsilon_{abc}J^c.
\eea
Here  indices are raised  with
the Lorentzian metric  $\eta^{ab}=(1,-1,-1)$ for Lorentzian signature but
with the trivial metric $\delta_{ab}$ for Euclidean signature.
To arrive at the Chern-Simons formulation we interpret the Cartan 
connection again as an Ehresmann connection on a bundle with
a bigger structure group.  
The structure groups which result for the various values of $\lambda$
and the two choices of signature are summarised in table 1. 

\vspace{1cm}

\noindent
\begin{tabular}{|c|c|c|}
\hline
  &    &   \\
Cosmological constant   & Euclidean signature & Minkowskian signature \\
  &    &   \\
\hline
  &    &   \\
$\lambda = 0$  & $ ISO(3)$ & $ISO(2,1)$ \\
 &    &   \\
\hline
  &    &   \\
$\lambda > 0$ & 
$ SO(4) \simeq  \frac{SU(2)\times SU(2)}{\ZZ_2}$&
$SO(3,1) \simeq SL(2,\CC)/\ZZ_2$ \\
  &    &   \\
\hline
  &    &   \\
$\lambda < 0$ & $ SO(3,1) \simeq SL(2,\CC)/\ZZ_2 $&
$SO(2,2) \simeq \frac{SL(2,\RR)\times SL(2,\RR)}{ \ZZ_2} $\\
  &    &   \\
\hline
\end{tabular}

\vspace{1cm}

\centerline{{\bf Table} 1}

\vspace{1cm}

Using the generalised commutation relations (\ref{comrels}) and 
again the bilinear form (\ref{inprod}) to 
expand the  Chern-Simons action (\ref{CSaction}) 
one indeed 
finds the  Einstein-Hilbert action for cosmological constant $\lambda$:
\bea
\label{EHcos}
S_{EH} [\omega,e,\lambda] = \int_M e_aF_\omega^a
+ \lambda \int_M \,\epsilon_{abc}\, e^a \wedge e^b \wedge e^c.
\eea

\section {Gravitational phase space and its symplectic  structure}

The phase space of a classical field theory is the space of 
solutions of the equations of motions, modulo gauge invariance.
Adopting the  Chern-Simons formulation of  three dimensional gravity 
we thus find that the phase space of gravity in  three  dimensions
is the moduli space of flat $G$-connection on the surface $\Sigma$,
where $G$ is the relevant gauge group extracted from table 1. 
In the following we specialise to compact Riemann surfaces of genus 
$g$. Also, we denote the  Lie algebra of $G$ by $\gothg$. 
Starting with the classic paper of Atiyah and Bott 
\cite{AB}, the  moduli space of flat $G$-bundles on $\Sigma$
 has been studied extensively for
  semi-simple, compact  Lie groups $G$. A pedagogical summary with
further references can be found in \cite{Atiyah}. 

The Lie groups relevant for
our discussion of gravity are typically not compact and only some are
semi-simple. This means in particular that the topology of various
quotient spaces  we are about to consider needs to be re-examined 
carefully. We shall restrict ourselves to 
a formal description of the relevant moduli spaces, but will highlight
some points which call for a more careful analysis. 
That  a formal  discussion is possible at all depends crucially on 
 the  existence of the 
non-degenerate invariant form (\ref{inprod}) on the Lie algebras
relevant for us.  Using that form, we briefly summarise the 
main results and formulae.

The tangent space to the  space of $G$ connections 
on a Riemann surface $\Sigma$ is the affine 
space of $\gothg$-valued one-forms on  $\Sigma$.
There is a natural symplectic form on this space, which can
also be derived from   the Chern-Simons action (\ref{CSaction}). 
It is 
\bea
\label{ABsympl}
\Omega = 
\int_{\Sigma}\langle \delta_1 A \wedge \delta_2 A \rangle .
\eea
One checks that it is invariant under the gauge transformations 
\bea
\label{gaugetrans}
A\rightarrow gAg^{-1} + dg g^{-1},
\eea
where $g\in $Map$(\Sigma,G)$.
An elementary calculation shows that 
the momentum mapping of this action is proportional to the curvature $F$.
Hence the space of  flat connections 
modulo gauge equivalence 
  is equal to the symplectic quotient of the space of all
 $G$-connections by the group of gauge transformations
(\ref{gaugetrans}).
This quotient 
is the moduli space of flat $G$-connections on $\Sigma$.
It  inherits a symplectic structure from the symplectic
structure on the space of all connections, which we refer to 
as the Atiyah-Bott symplectic structure. 
There are various ways of describing
the moduli  space and its symplectic structure. 
For us, the following description in terms of representations 
of the fundamental group is   most useful. One of its advantages
is that it is straightforward  to include marked points on the 
Riemann surface. In the following we thus write 
$\Sigma_{g,m}$ for a compact Riemann surfact of genus $g$ with
$m$ marked points.

Let $\pi$ be the fundamental group of $\Sigma_{g,m}$. 
This group is generated by $2g+m$
invertible
generators $a_1,b_1,...a_g,b_g,l_1,...l_m$ satisfying the relation
\bea 
[b_g,a_g^{-1}] ...[b_1,a_1^{-1}]l_m...l_1 =1,
\eea
where $[x,y]=xyx^{-1}y^{-1}$.
A flat $G$-connection on $\Sigma_{g,m}$
associates to each generator a holonomy element in $G$. However,
the insertion of the  charges 
at marked points as in (\ref{gravwith})
means that the holonomy around the i-th marked point is forced 
to lie in a fixed conjugacy class ${\mathcal C}_i$.
We shall discuss in detail how the conjugacy classes are related
to the co-adjoint orbits of the $\gothg^*$-elements  (\ref{dualphase})
in the next section.

The 
moduli space of flat connections on $\Sigma_{g,m}$
 depends on these conjugacy classes,
as well as on the group $G$ and  the genus $g$. We define it as 
\bea
M(G,g,\calC_1,...\calC_m) =\{\rho\in  \mbox{Hom}(\pi,G), \rho(l_i)\in 
\calC_i\}/G,
\eea
where  the group $G$ acts via conjugation.
In the case where $G$ is compact and semi-simple this space is 
a Hausdorff space, and it is a manifold at irreducible points
(where the image of $\pi$ generates $G$). 
It would clearly be important to know if similar results 
hold for the groups $G$ in table 1. 

For our approach to the eventual  quantisation of this space
it is important that it can also be written as a quotient
of the space  $G^{2g}\times \calC_1\times ...\calC_m$. 
To see this we introduce the group-valued momentum
mapping  $\mu:G^{2g}\times \calC_1\times ...\calC_m \rightarrow G$,
\bea
\mu(A_1,B_1, ....A_g,B_g,L_1,...L_m)= 
[B_g,A_g^{-1}] ...[B_1,A_1^{-1}]L_m...L_1.
\eea
Then 
\bea
M(G,g,\calC_1,...\calC_m)=\mu^{-1}(1)/G.
\eea

The important  observation of 
of Fock and Rosly \cite{FR} is that, under certain conditions,
the space $G^{2g +m}$ has  a natural Poisson structure. When restricted to 
conjugation invariant functions on $G^{2g}\times \calC_1\times ...\calC_m$
(which therefore descend to functions 
on $M(G,g,\calC_1,...\calC_m)$) that Poisson structure 
  agrees with the Poisson structure derived 
from the Atiyah-Bott symplectic structure on  the moduli space
$M(G,g,\calC_1,...\calC_m)$. To write down the Fock-Rosly  Poisson structure
one requires an element $r\in\gothg\otimes\gothg$
which satisfies  the classical Yang-Baxter equation and is such
that its symmetric part agrees with the non-degenerate invariant form 
$\langle\,\,,\,\,\rangle$ used
in the definition of the symplectic structure (\ref{ABsympl}).
 In the next section we shall write down such an $r$-matrix 
for  $\gothg=iso(3)$. 

\section{The Lie bi-algebra structure of $iso(3)$}

In the following discussion 
 we use  the conventions and terminology of 
\cite{CP}. In particular we 
 write $\sigma:\gothg\otimes\gothg\rightarrow \gothg\otimes\gothg$
for the flip  operation
 $\sigma (X\otimes Y) =Y\otimes X$, and 
if $r=r^{\alpha\beta}
X_\alpha \otimes X_\beta\in \gothg\otimes\gothg$, then  $r_{12}\in
\gothg\otimes\gothg\otimes\gothg$ is defined to  be 
$r_{12}=r^{\alpha\beta} X_\alpha \otimes X_\beta\otimes 1$.

 Recall that a
 bi-algebra structure on a Lie algebra $\gothg$ is a skew-symmetric
 linear map $\delta:\gothg  \rightarrow \gothg \otimes \gothg$,
called the co-commutator, which satisfies the co-cycle condition
$$
\delta[X,Y]= (\mbox{ad}_X\otimes 1 + 1\otimes\mbox{ad}_X)
 \delta(Y))-(\mbox{ad}_Y\otimes 1 + 1\otimes\mbox{ad}_Y) \delta(X)
$$ 
and is such  that 
the dual  $\delta:\gothg^* \otimes \gothg^*\rightarrow \gothg^*$
is a commutator. 
A Lie bi-algebra $\gothg$  is  co-boundary if the 
co-commutator $\delta$ can be written as 
$\delta(X) =(\mbox{ad}_X\otimes 1 + 1\otimes\mbox{ad}_X)  r$,
for some element  $r\in \gothg\otimes \gothg$ which satisfies
two  conditions.
The first is that $r+\sigma(r)$
is an invariant element of $\gothg\otimes \gothg $. 
 The second condition is that 
\bea
[[r,r]]:=[r_{12},r_{13}]+[r_{12},r_{23}]+[r_{13},r_{23}]
\eea
is an invariant element of $\gothg\otimes \gothg\otimes \gothg$
When that invariant element is $0$, $r$ is said to satisfy 
the classical Yang-Baxter equation (CYBE) and the bi-algebra is 
said to be quasi-triangular.

To  exhibit the relevant Lie bi-algebra structure of $iso(3)$ we first
recall   the commutation relations for the rotation generators $J_a$ and 
the translation generators $P_a$ of  $iso(3)$ given in (\ref{spinrels})
and (\ref{transrels}). Further, we identify 
the vector space dual $iso(3)^*$ with $iso(3)$ via 
the non-degenerate pairing (\ref{inprod}) on $iso(3)$.
 Note that the basis elements dual to the rotation generators
are the translation generators and {\it vice-versa}.

The Lie algebra  $\gothg =
iso(3)$ has several bi-algebra structures \cite{Stachura}.
The one which is relevant for us is associated to the interpretation of 
$iso(3)$  as the classical double of the Lie algebra $so(3)$ (see \cite{CP}
for a discussion of classical doubles).
Note that this double structure is different from the  bi-algebra structure
one obtains by thinking of $iso(3)$ as a Wigner contraction of $so(4)=
so(3)\oplus so(3)$ and using the standard bialgebra structure of $so(3)\simeq
sl_2$. Quantising the latter leads to the q-deformed universal enevloping
algebra $U_q(iso(3))$ as described for example 
in \cite{CG}. Quantising $iso(3)$
viewed as a classical double yields, as we shall see, a different 
quantum group, namely the quantum double
$D(SU(2))$. For a classical double, there is a particularly simple
formula for  the classical $r$-matrix. In our case it is  
\bea
\label{rmatrix}
r= P_a\otimes J_a.
\eea
One easily checks that it satisfies the CYBE and that 
the symmetrised part 
\bea 
r^s:=\frac{1}{2}(r+ \sigma(r)) = \frac{1}{2}( P_a\otimes J_a
+ J_a\otimes P_a)
\eea
 is invariant under the adjoint action of
$iso(3)$. Note in particular that  $r^s$ corresponds to
the non-degenerate pairing (\ref{inprod}) on the vector 
space $iso(3)$. We conclude that the bi-algebra  $iso(3)$ is co-boundary
and quasi-triangular, and that its  co-commutator  $\delta: iso(3)\rightarrow iso(3)\otimes iso(3)$ is 
\bea
\label{delmap}
\delta(P_a)=\epsilon_{abc} P_b\otimes P_c, \qquad \delta (J_a)=0.
\eea
Using the pairing (\ref{inprod}) we find that the dual map $\delta^*:
 iso(3)^*\otimes iso(3)^*\rightarrow  iso(3)^*$ is   
\bea
\label{coiso}
\delta^*(J_a,J_b)=\epsilon_{abc} J_c \qquad
\delta^*(J_a,P_b)= 
\delta^*(P_a,P_b)=0.
\eea
This defines the Lie algebra structure of $iso(3)^*$.
We deduce the Lie-algebra isomorphism  $iso(3)^*\simeq so(3)\oplus \RR^3$.

The classical $r$-matrix (\ref{rmatrix}) is the only
input needed to write down  the Fock-Rosly  Poisson structure on 
the phase space of three-dimensional Euclidean gravity.
To review the Fock-Rosly structure briefly in a more general setting,
 let $G$ be a 
 simply connected Lie group with Lie algebra  $\gothg$.
Fock and Rosly first gave a general formula for the Poisson bracket
of functions on $G^{2g+m}$ in \cite{FR}.
A  particularly convenient form can be found in \cite{AS}. It 
involves the vector fields generated by left- and right action of the 
group $G$ on $G^{2g+m}$. 
The general expresssion  is quite complicated and we will not write 
it down here.  It is, however,  instructive
to write it down
in the simplest case $g=0$ and $m=1$. In that case we obtain 
a Poisson structure on the group itself. 
Let $\{X_\alpha\}_{\alpha=1,...\mbox{dim}\gothg}$ be a basis of $\gothg$
and write $ X^L_\alpha$ and $X^R_\alpha$ for the vector fields on $G$
generated by left- and right action of $\gothg$ on $G$. For a function
$f$ on $G$ we define more precisely
\bea
X^L_\alpha f(g)=\frac{d}{dt}f (e^{-tX_\alpha}g)|_{t=0} \quad
\mbox{and} \quad
X^R_\alpha f(g)=\frac{d}{dt}f (ge^{tX_\alpha})|_{t=0}.
\eea
Then, if $r=r^{\alpha\beta}X_\alpha\otimes X_\beta$ is a solution 
of the CYBE,  the Fock-Rosly bracket of two functions $f_1$ and $f_2$ 
on $g$ is 
\bea
\label{frbr}
\{f_1,f_2\}&=& \frac{1}{2}r^{\alpha\beta}\bigl(X^R_\alpha f_1 X^R_\beta f_2
+X^L_\alpha f_1 X^L_\beta f_2 - 
X^R_\beta f_1 X^R_\alpha f_2 - X^L_\beta f_1 X^L_\alpha f_2\bigr)
\nonumber \\
&+&r^{\alpha\beta}\bigl(X^R_\alpha f_1 X^L_\beta f_2-
X^L_\beta f_1X^R_\alpha f_2\bigl).
\eea

To compute this bracket for the case   $\gothg=iso(3)$, 
we  need some notation.  The double cover  $G=ISO(3)^\sim$
is the semi-direct product $\RR^3\rtimes SU(2)$.  We write elements
of $ISO(3)^\sim$ as pairs $(\ba,u)$ 
 of translation vectors  $\ba\in \RR^3$ and 
elements $u\in SU(2)$.
For an element $u\in SU(2)$ we write $\ad (u)$ for the $SO(3)$ matrix 
representing $u$  in the adjoint representation.  Then the multiplication law
for elements $(\ba,u),(\bob,v)\in ISO(3)^\sim$ is
\bea
\label{isogroupmult}
 (\ba,u)\cd (\bob,v)=(\ba +\ad(u)\bob, uv).
\eea
In the present context it is useful to 
parametrise group elemens via the exponential map.
More precisely, we  first identify 
an element $\xi^* \in iso(3)^*$  with  an element $\xi \in iso(3)$
via the inner product (\ref{inprod}) and then map it into $ISO(3)^\sim$
via the exponential map.
The combined map is 
\bea
\label{funexp}
\exp\circ*&:& iso(3)^* \rightarrow ISO(3)^\sim \nonumber \\
&& p_a P^*_a +j_a J^*_a \mapsto \exp(j_a P_a +p_a J_a).
\eea 
To make contact with our previous parametrisation of 
$ISO(3)^\sim$ elements  in terms of translations $\ba$ and 
rotations $u$ note that if $\exp(j_a P_a +p_a J_a)$ is the element 
$(\ba,u)$, then
$u=\exp(p_aJ_a)$ and $\ad(u^{-1})\ba=\bj$. Since the exponential
map  from $so(3)$ to $SU(2)$ is onto, it follows 
that the exponential map (\ref{funexp}) is onto.  Recalling our 
description of the co-adjoint orbits of $ISO(3)^\sim$ following eq. 
(\ref{phaseco}) we conclude that the orbits ${\mathcal O}_{\mu s}$ with 
$0\leq\mu<2\pi$ are mapped bijectively to conjugacy classes 
in $ISO(3)\sim$, 
which we denote by $\calC_{\mu  s}$. The only conjugacy classes
in $ISO(3)\sim$ which are  not bijectively related to  co-adjoint orbits
in this way are  the classes containing  elements  with $u=-1$.

The vector fields $J^L_a$ and $P^L_a$ generated by 
left-action of, respectively,  the rotations and translations  
on $ISO(3)^\sim$ are  given by 
\bea 
J^L_a f (g) &=& \frac{d}{dt}f (e^{-tJ_a}g)|_{t=0} \nonumber \\
&=& \bigl(\frac{\ad(u)}{1-\ad(u)}\bigr)_{ab} \epsilon_{bcd} p_c
\frac{\partial f }{\partial p_d}(g)
\eea
and 
\bea 
P^L_af (g)& =& \frac{d}{dt}f (e^{-tP_a}g)|_{t=0}\nonumber \\
&=& -\ad(u)_{ab}\frac{\partial f } {\partial j_b}(g).
\eea
In this and  the following formulae, 
all functions are evaluated at  
 $g=\exp(p_a J_a +j_a P_a)$.
The vector fields generated by right actions of rotations and translations 
are 
\bea
J^R_a f (g)& = &\frac{d}{dt}f (ge^{tJ_a})|_{t=0} \nonumber \\
&=&\bigl( \frac {1} {\ad(u)-1}\bigr)_{ab}\epsilon_{bcd}p_c \frac{\partial f }
{\partial p_d}(g) -\epsilon_{abc}j_b\frac{\partial f }{\partial j_c}(g)
\eea
and
\bea 
P^R_af (g) = \frac{d}{dt}f (ge^{tP_a})|_{t=0}
=\frac{\partial f }{\partial j_a}(g).
\eea

The Fock-Rosly Poisson bracket (\ref{frbr}) of functions $f_1$ and $f_2$ on 
$ISO(3)^\sim$ now takes 
the following form.
\bea
\{f_1,f_2\}& =&
\frac{1}{2}\left( P^R_a f_1 J^R_a f_2 +  P^L_a f_1 J^L_a f_2 -
 J^R_a f_1 P^R_a f_2-  J^L_a f_1 P^L_a f_2\right) \nonumber \\
&&+P^R_a f_1J^L_a f_2-
J^L_a f_1P^R_a f_2 \nonumber \\
&=&\epsilon_{abc}\left(
\frac{\partial f_1}{\partial p_a}\frac{\partial f_2}{\partial j_b} p_c
+\frac{\partial f_1}{\partial j_a}\frac{\partial f_2}{\partial p_b} p_c \right)
+\epsilon_{abc}\frac{\partial f_1}{\partial j_a}\frac {\partial f_1}
{\partial j_b} j_c.
\eea
In particular we find  that the coordinate functions 
$j_a$ and $p_a$, $a=1,2,3$,  have the Poisson brackets (\ref{kiri}).
 We therefore
conclude that the restriction of the Fock-Rosly Poisson bracket
to a conjugacy classes 
$\calC_{\mu s}$ in $ISO(3)^\sim$ with  $\mu\not \in 2\pi \ZZ$
is the push-forward of the 
Kostant-Kirillov Poisson  structure on ${\mathcal O}_{\mu s}$ via the map 
(\ref{funexp}).

\section{The quantum double $D(SU(2))$ as a deformation of 
the group algebra of $ISO(3)^\sim$}

The quantum double of a finite group $H$ is a quasi-triangular Hopf
algebra constructed, via Drinfeld's  double construction,
out of the Hopf algebra of functions on $H$. As a vector space
$D(H)$ is the tensor product of the  algebra of functions $F(H)$
on $H$ and the group algebra $\CC(H)$. The quantum double
$D(H)$  plays an important role in the physics of 
orbifold conformal field theories \cite{DVVV} and in discrete gauge theories 
\cite{BDP}. Its mathematical structure as a quasi-triangular Hopf
algebra was first discussed in \cite{DPR}. When  generalising the 
construction  to  locally compact Lie groups $H$,  one 
has to choose appropriate sets of (generalised) functions on $H$ and 
different choices have been made by different authors.
 For a construction which emphasises
the duality between $D(H)$ and $D(H)^*$  see \cite{Bonneau} and for
an approach via von Neumann algebras see \cite{Mueger}. 
In \cite{KM} a definition of $D(H)$ is given which emphasises
the similarity with transformation group algebras. That  approach 
is particularly well-adapted to the classification 
of the irreducible representations of $D(H)$, see \cite{KBM}.
For our purpose it is best to 
use a definition which is
 close to that of \cite{KM} but slightly more general. 
In this paper we are exclusively concerned with the case  $H=SU(2)$. 
The main
reference is \cite{BM}, where the role of   $D(SU(2)$ as the underlying
symmetry of  quantised $ISO(3)^\sim$ Chern-Simons theory was deduced from
a detailed study of the quantum numbers, fusion and braiding properties of 
particles in that theory. However, the question of how to obtain the 
quantum group $D(SU(2))$ from a 
quantisation of the Poisson structure on the phase
 space  of $ISO(3)^\sim$ Chern-Simons theory   
was not clarified. The goal of this section
is to fill that gap. The chief difficulty which we have to overcome
is the identification of the relevant deformation parameter. As we
shall see presently,  this is not manifest in the formulation of $D(SU(2))$.

We begin with a few general definitions. For a locally compact, unimodular 
group $H$ with Haar measure $dh$ we define   $M(H)$ to be 
set of complex-valued,  bounded  measures on $H$ which are 
 absolutely continuous   with respect to $dh$ or
 pure point  measures (thus we exclude continuous singular measures).
For every 
absolutely continuous measure $\mu$ on $H$ there is a function
$f\in L^1(H,dh)$ such that  $d\mu =f dh$. Pure point measures
have no such representation, but for our calculations it is convenient
to adopt the physicists' habit of  writing pure point measures
as a product of a Dirac delta function and the Haar measure.
With that convention, the measures in $M(H)$ have the decomposition
\bea
\label{meadecomp}
d\mu = (f  + \sum_{i\in \NN}\lambda_i  \delta_{h_i})dh,
\eea
where  $f\in L^1(H,dh)$, $\{\lambda_i\}_{i\in \NN}\in l^1(\CC)$
and $\delta_{h_i}$ is the normlised Dirac delta function at $h_i$.
$M(H)$ is a Banach algebra with norm
\bea
||\mu||_1 = \int_H |d\mu|
\eea 
and multiplication defined by convolution. The convolution 
can be defined without reference to the decomposition (\ref{meadecomp})
(see e.g. \cite{Pedersen}), but using it simplifies the notation.
Representing $d\mu_1$ and $d\mu_2$ 
 by   generalised functions $f_1$, and $f_2$ which 
are   in $L^1(H)$, or   delta functions  or a linear combination
of the two, the convolution becomes  
the ordinary convolution of functions:
\bea
f_1*f_2(h)=\int_H dv f_1(v)f_2(v^{-1}h).
\eea

The algebra $M(H)$ will play the role of the group algebra  of 
$H$ in the rest of the paper. The role of the  function algebra
will be played by bounded, uniformly continuous functions $C_B(H)$ on $H$.
When equipped with the supremum norm this, too, is  a Banach algebra 
Then we define the  quantum double $D(H)$ of a locally compact
 Lie group  $H$ with Haar measure $dh$
 as the set of  bounded, $C_B(H)$-valued  measures on $H$ 
 which are either  absolutely
continuous with respect to $dh$ or pure point measures. 
More practically we think of 
elements of $D(H)$ as  generalised functions $F$ on $H\times H$ which 
are bounded and continuous in the first argument and 
  measures of the form (\ref{meadecomp})  in the second argument.
The norm of such a function is 
 \bea
||F||_1=\int_H dv \,\,\mbox{sup}_{h\in H}|F(h,v)|.
\eea
We refer the reader to
\cite{KM} and \cite{BM} for a complete list of how  the algebraic operations
of a quasi-triangular Hopf-$*$-algebra are implemented in $D(H)$.
Note that the class of functions considered  there is smaller than
the ones included in our definition of $D(H)$. Some of the formula
in \cite{KM} involve $\delta$-functions and only make rigorous 
sense in our definition of $D(H)$.
For our purposes we only need the formulae  for the multiplication
and the co-multiplication. The multiplication of two elements
$F_1$ and $F_2$ of $ D(H\times H)$ is 
\bea
\label{multi}
\bigl(F_1\bullet F_2\bigr)(h,u)=\int_H\, dw\, F_1(h,w)F_2(w^{-1}hw,w^{-1}u).
\eea
The co-multiplication $\Delta$  is defined via 
\bea
\label{comulti}
\Delta F(h_1,u_1,h_2,u_2) = F(h_1h_2,u_1) \delta_{u_1}(u_2).
\eea
The final ingredient we need is the universal R-element:
\bea
\label{bigrmatrix}
R(h_1,u_1,h_2,u_2)=\delta_{h_1}(u_2)\delta_e(u_1).
\eea

The irreducible representations (irreps) 
 of $D(SU(2))$ are classified in \cite{KM}  and
 are labelled by pairs $(\mu,s)$, where $\mu\in[0,2\pi]$  labels 
a conjugacy class $C_\mu$ in $SU(2)$  and $s$ labels an irrep of the 
centraliser of a chosen point on the conjugacy class. Explicitly,
we realise the rotation generators $J_a$ as $J_a=-\frac{i}{2}\tau_a$,
where the $\tau_a$ are the Pauli matrices. Then every element $h\in
 SU(2)$ can be written in axis-angle parametrisation as 
\bea 
h(\mu,\bn) = \exp(\mu (n_aJ_a))
\eea 
for a unit vector $(n_1,n_2,n_3)$ and $\mu\in [0,2\pi]$.
There are two single-element conjugacy classes,
$C_1=\{1\}$  and $C_{-1}=\{-1\}$. In both cases 
the centraliser group is the whole of $SU(2)$, whose irreps  are labelled
by a positive half-integer, the spin.  The other conjugacy
classes are topological two-spheres, 
consisting of all rotations by some angle $\mu\in(0,2\pi)$ about an 
arbitrary axis: $C_\mu=\{h(\mu,\bn)|\bn\in S^2\}$.
In each class  $C_\mu$  we single out   
the rotation about the 3-axis 
$h_\mu:=\exp(\mu J_3)$   to define  the centraliser group belonging to
$C_\mu$. That group therefore consists of all rotations around the
3-axis  and is isomorphic to $U(1)$. Its representations are 
labelled by a half-integer, which one may think of as a ``Euclidean
 helicity'' in the present context. We denote them by $\pi_s$.
For $\mu\in(0,2\pi)$ and $s\in\ZZ/2$  let
\bea
\label{Hilbert}
H_{s}=\{ \phi\in L^2(SU(2),\CC)|\phi(x h_{\omega})= \pi_s(h_\omega^{-1})
\phi(x) \}.
\eea
This is  the carrier space of the representation $\Pi_{\mu s}$. The
 action of an element $F\in D(SU(2))$ is
\bea
(\Pi_{\mu s}(F)\phi) (x) = \int_{SU(2)}\,dw\,
 F(x h_\mu x^{-1},w)\phi(w^{-1}x).
\eea
The irreps labelled by $(0,j)$, $j=0,1/2,1 ...$ are the standard  spin $j$
 representations of $SU(2)$ with carrier space $\CC^{2j+1}$.
In order to simplify formulae in the following sections it is 
convenient to realise this carrier space as a certain space $H_{j}$   of 
functions on $SU(2)$. Writing $D^j(x)$ for the matrix
representing  $x\in SU(2)$ in the standard 
spin $j$ representation and 
 $|j,j\rangle$ for the highest weight vector
in that  representation, we define a function $\phi$ 
for every $\varphi \in \CC^{2j+1}$ via 
$\phi(x)=\langle j,j|D^j(x^{-1})|\varphi\rangle$. $H_{j}$ is the linear 
space of functions obtained from $\CC^{2j+1}$  in this way. Equivalently
one could define it as the span of the  
 Wigner functions $D_{nj}^j$,  $n=-j,-j+1,...,j-1,j$ on
 $SU(2)$. 
Then the action of an element $F\in D(SU(2))$ on $\phi\in H_{j}$  is
\bea
(\Pi_j(F)\phi)(x) = \int_{SU(2)}\,dw\,
 F(1,w)\phi(w^{-1}x).
\eea
Finally we note  the action of the universal $R$-element
(\ref{bigrmatrix}) in the tensor product representation 
 $\Pi_{\mu_1 s_1}\otimes\Pi_{\mu_2 s_2}$   on  some state 
$\Phi \in   H_{s_1} \otimes H_{s_2}$. It 
is 
\bea
\label{Raction}
\bigl( \Pi_{\mu_1,s_1}\otimes \Pi_{\mu_2,s_2} \,(R)\,\Phi \bigr)(x_1,x_2))
= \Phi(x_1,x_1h_{\mu_1}x_1^{-1}x_2).
\eea

In order to establish the 
relationship between $D(SU(2))$and $ISO(3)^\sim$  we
need to  review some properties of $ISO(3)^\sim$. As a manifold
 $ISO(3)^\sim\simeq \RR^3\times SU(2)$, and the Haar measure
on $ISO(3)^\sim$ is the product of the Lebesgue measure on $\RR^3$
and the Haar measure on $SU(2)$. Concretely, if $g=(\ba,u)\in ISO(3)^\sim$
then we use the Haar measure $dg=d^3\ba\,du$. Then  we can realise
the group algebra of $ISO(3)^\sim$ as described earlier  for general
locally compact groups as the set $M(ISO(3)^\sim)$ of bounded
 measures on $ISO(3)$ either absolutely continuous with
respect to $dg$ or pure point.    Again representing such
measures by generalised functions  $\hat f_1,\hat f_2$
(possibly including delta-functions)
the multiplication rule
is the convolution
\bea
(\hat f_1\hat \bullet\hat  f_2) (\ba,u) =\frac{1}{2\pi^3}
\int_{\RR^3\times SU(2)} 
d^3b\,dw\,\hat f_1(\bob,w)\hat f_2(\ad(w^{-1})(\ba-\bob),w^{-1}u).
\eea
Now we perform  a  Fourier transform on the first argument of $\hat f$,
thus obtaining a function $f$ on $(\RR^3)^*\times SU(2)$:
\bea
\label{halff}
 f(\bk,x)=\frac{1}{(2\pi)^3}\int_{\RR^3}d^3a\,
 \exp({-i\bk\cd\ba})\hat f(\ba,x).
\eea
Elements of the dual space  $(\RR^3)^*$ have the dimension of 
inverse length. Eventually they will physically be interpreted  
as momenta. Strictly speaking this is only possible after
the introduction of a constant of the dimension length$\times$momentum,
such as Planck's constant. Anticipating a more detailed discussion of 
dimensions and physical interpretations in the next section 
we refer to $(\RR^3)^*$ as momentum space from now onwards.

The Fourier transform of a bounded measure on $\RR^3$  is 
a bounded,  uniformly continuous function on 
 $(\RR^3)^*$.
  Although not all bounded,  uniformly continuous functions  
are obtained in this way, we define the group algebra of $ISO(3)^\sim$
as the set of of bounded measures on $SU(2)$, absolutely continuous
 with respect
to the Haar measure or pure point,
taking values in the Banach space of bounded,  
uniformly continuous functions on $(\RR^3)^*$. 
As in the case of $D(SU(2))$
we think concretely of generalised 
functions  $f$ on $(\RR^3)^*\times SU(2)$ which are bounded,  
uniformly continuous functions of the first argument and 
measures of the form (\ref{meadecomp}) with respect to the second.
For reasons which will become clear later we denote the set
of all such generalised functions by $A_0$. We give $A_0$ 
the structure of a Hopf algebra as follows.  The product 
of two generalised 
functions $f_1,f_2$ in $A_0$ is obtained by applying the 
 Fourier transform (\ref{halff}) to the convolution
product. The result is
\bea
\label{isomult}
\bigl(f_1\bullet f_2\bigr) (\bk,u)=\int_{SU(2)}\,dw\, f_1(\bk,w)
f_2(\ad(w^{-1})\bk,w^{-1}u).
\eea
The group-like co-multiplication for $ISO(3)^\sim$
leads to  the following  co-multiplication  for $f\in A_0$ 
\bea
\label{isocomult}
(\Delta f)(\bk_1,u_1,\bk_2,u_2)=f(\bk_1+\bk_2,u_1)\delta_{u_1}(u_2).
\eea
One checks that with antipode
\bea 
\label{antipode}
S f(\bk,u)=f(-\ad(u^{-1})\bk,u^{-1}), 
\eea
unit
\bea
\label{unit}
1(\bk,u)=\delta_e(u)
\eea
and co-unit
\bea
\label{co-unit}
\epsilon(f)=\int_{SU(2)}f(0,u)
\eea
$A_0$ is a co-commutative Hopf algebra.

The irreps of the Euclidean group $ISO(3)^\sim$ and hence also of
$A_0$  are labelled 
by pairs of $SU(2)$ orbits and centraliser representations, much 
like for the double  $D(SU(2))$.
Now the relevant orbits are the orbits of the $SU(2)$ action on
$(\RR^3)^*$. There is one orbit consisting of the origin $(0,0,0)$,
and all the other orbits are two-spheres in $(\RR^3)^*$. These orbits 
are therefore 
naturally  labelled by their   radius.
In the current context this radius has the interpretation of Euclidean
mass, so we denote  it by  $M\in\RR^{\geq 0}$. If $M=0$ the 
centraliser group is the whole of $SU(2)$, with irreps labelled by
 the spin $j$ as before.  For $M>0$ we single out the point
$\bd^M_3=(0,0,M)$ on the $3$-axis. Then the centraliser group is 
the group of rotations around the $3$-axis. The generic irreps
of   $ISO(3)^\sim$ are therefore labelled by pairs $(M,s)$, 
with $M>0$ and $s\in \ZZ/2$ labelling an irrep of $U(1)$.
The carrier spaces of  such irreps  are the Hilbert spaces 
$H_s$  defined in (\ref{Hilbert}). The action of $f\in A_0$ is
\bea
\label{isorep}
\bigl(\Pi_{Ms}(f)\phi\bigr)(x) =
\int_{SU(2)} \, dw \,f(\bk(x,M),w)\phi(w^{-1}x),
\eea
where 
\bea
\label{bpdef}
\bk(x,M)=\ad(x)\bd^M_3.
\eea 
Note that elements  $(\ba,u)\in ISO(3)^\sim$ correspond to 
 plane waves in momentum space and $\delta$-functions
in $SU(2)$:
$f_{\ba,u}(\bk,v))=\exp(i\ba\cd\bk)\,\delta_u(v) $.
Applying the formula (\ref{isorep}) to such functions 
one obtains the familiar formula for the representation of 
$(\ba,u)\in ISO(3)^\sim$
\bea 
\bigl(\Pi_{Ms}((\ba,u))\phi\bigr)(x) 
= \exp(i\ba\cd\bk(x,M))\phi(u^{-1}x).
\eea
 
From the representation  of 
the Lie group, we obtain representations of the Lie algebra
in the usual way.  Writing $\bd_a$, $a=1,2,3$,  for 
the canonical basis of $(\RR^3)^*$ ($\bd_1=(1,0,0)$ etc. ) 
we have the following representation of  generators $P_a$
of translations:
\bea
\label{transrep}
\bigl(\Pi_{Ms}(P_a)\phi\bigr)(x)=
\frac{d} {d \epsilon}
\bigl(\Pi_{Ms}((\epsilon\bd_a,1))\phi\bigr)(x)|_{\epsilon=0}=
i \,k_a(x,M) \phi(x).
\eea
For the rotation generators $J_a$ we find
\bea
\bigl( \Pi_{Ms}(J_a)\phi\bigr)(x)=
\frac{d} {d \epsilon}
\bigl(\Pi_{Ms}((0,\exp(\epsilon J_a)))\phi\bigr)(x)|_{\epsilon=0}, 
= -{\mathcal J}_a \phi(x)
\eea
where ${\mathcal J}_a$, $a=1,2,3$,  are the vector fields on $SU(2)$  
which generate the 
left action   of  $SU(2)$  on itself.
We deduce in particular that the  classical $r$-matrix (\ref{rmatrix})
acts on elements $\Phi$ of the  tensor product  $H_{s_1}\otimes H_{s_2}$ via
\bea
\label{raction}
\bigl( (\Pi_{M_1s_1}\otimes \Pi_{M_2s_2})(r) \Phi\bigr)(x_1,x_2)
= -i\bigl( k_a(x_1,M_1){\mathcal J}^{(2)}_a\bigr) \Phi(x_1,x_2),
\eea
where we use the superscript $^{(2)}$ to indicate that ${\mathcal J}^{(2)}$
acts on the second argument of $\Phi$.

The final ingredient we need in order  establish the   promised 
relationship between $D(SU(2))$ and $ISO(3)^\sim$ is a family
of exponential maps from the Lie algebra $so(3)$
to $SU(2)$. Let $J^\kappa_a=\kappa J_a$  for a real, positive
parameter  $\kappa$, so that 
\bea 
[J^\kappa_a,J^\kappa_b] = \kappa
\epsilon_{abc}J^\kappa_a.
\eea
Then define
\bea
\label{expfun}
\exp_\kappa:
(\RR^3)^*\rightarrow SU(2),\quad \exp_{\kappa}(\bk)=\exp(k_aJ^\kappa_a).
\eea
This map is injective when restricted to the open Ball $B^3_{(2\pi/ \kappa)}$
 of radius $2\pi/ \kappa$ in $(\RR^3)^*$, but it is not surjective:
the image of $B^3_{(2\pi /\kappa)}$
is $SU(2)\backslash \{-1\}$. 
Now let $\hat \bk=\bk/|\bk|$ for $\bk\neq 0$, and write
$S^2_{(2\pi/\kappa)}$
for the two-sphere of radius $2\pi/\kappa$ in $(\RR^3)^*$.  Then
 define the space 
\bea
V_\kappa= \{f\in C_B((\RR^3)^*)|f|_{S^2_{2\pi/ \kappa}}=
\mbox{const},\, f(\bk- (4\pi/\kappa) \hat \bk)=f(\bk)   
\,\forall 
\bk\in\RR^3\backslash \{0\} \}.
\eea
Note that any $f\in V_\kappa$ necessarily takes the same value on
all two-spheres of radius $4m \pi/  \kappa  $, $m\in\ZZ$ as at the origin.
For us $V_\kappa$ is important because it contains the image  of the 
pull-back map
\bea
\label{exmap}
\exp_\kappa^*:C_B(SU(2))\rightarrow C_B((\RR^3)^*).
\eea
This map is injective and its inverse can be defined if 
we restrict it to $V_\kappa$. We denote the inverse by 
\bea
\label{lomap}
\log_\kappa^*:V_\kappa\rightarrow C_B(SU(2)).
\eea
Finally we define, for fixed $\kappa$, the space 
$A_\kappa$ as the space of  $V_\kappa$-valued,
  bounded, absolutely continuous or pure point
measures on $SU(2)$.  Again we
think of elements of $A_\kappa$ concretely as functions on $(\RR^3)^*
\times SU(2)$ which are in $V_\kappa$ as function of the first
argument and of the form (\ref{meadecomp}) in the second argument.  
It is clear that for each positive 
value of $\kappa$, $A_{\kappa}\subset A_0$ and that $A_\kappa$
inherits an algebra structure from $A_0$. In fact $A_\kappa$ 
is a consistent truncation  of $A_0$  as  an  algebra:
$A_{\kappa}$ is closed under the multiplication $\bullet$
(\ref{isomult}) and the unit (\ref{unit}) of $A_0$  
is contained in $A_\kappa$. However, 
the co-multiplication $\Delta$ (\ref{isocomult})
does not map $A_{\kappa}$ into  $A_{\kappa}\otimes A _{\kappa}$, so
that  $A_{\kappa}$ is not  a consistent truncation 
of $A_0$  as  a Hopf   algebra. We shall now 
 show how to turn $A_\kappa$ into a quasi-triangular Hopf algebra
by introducing a deformed  co-multiplication.

To construct the new co-multiplication
we  extend the map (\ref{exmap}) to a map
\bea
\mbox{EXP}_\kappa^*: D(SU(2))\rightarrow A_\kappa,
\eea
by pull-back on the first argument.
This map is again invertible, and we denote the inverse by
\bea
\mbox{LOG}_\kappa^*:A_\kappa\rightarrow D(SU(2)).
\eea
Now we can use the $\kappa$-dependent bijection EXP$_\kappa$ to induce
a Hopf-algebra structure from $D(SU(2))$
on $A_\kappa$.
This re-produces the $\kappa$-independent   multiplication  rule
 (\ref{isomult}), the antipode (\ref{antipode}), unit (\ref{unit}) and 
co-unit (\ref{co-unit}).
The co-multiplication $A_\kappa \rightarrow A_\kappa
\otimes A_\kappa$, however,  does depend on $\kappa$:
\bea
\label{newco}
(\Delta_\kappa f)(\bk_1,u_1,\bk_2,u_2)= \mbox{LOG}_\kappa^*\bigl(
f(\exp_\kappa(\bk_1)\exp_\kappa(\bk_2),u_1)\bigr)\,\,
\delta_e(u_1^{-1}u_2).
\eea

For a better understanding of this formula and of 
the limit we will consider shortly, it is useful to think of it 
 in terms of the  
Baker-Campell-Hausdorff expression  for the product of exponentials
of Lie algebra elements $X,Y\in su(2)$:
$\exp( X) \exp( Y)=\exp (X*Y)$, where
\bea 
X*Y=X+Y+\frac{1}{2}[X,Y]+ \frac{1}{12}[X,[X,Y]] + 
\frac{1}{12}[[X,Y],Y] +... .
\eea
With $X=k^a_1J^\kappa_a$,$ Y=k_2^a J^\kappa_a$, this induces a
$\kappa$-dependent
 product on $(\RR^3)^*$:
\bea
\label{BCH}
\bk_1*\bk_2= \bk_1+\bk_2 +\frac{\kappa}{2}\bk_1\wedge\bk_2 +\frac{\kappa^2}
{12 }
\bigl( \bk_1\cd\bk_2(\bk_1+\bk_2)-\bk_1^2\bk_2 -\bk_2^2\bk_1\bigr)
+... .
\eea
This is the product used in the first argument of $f\in A_{\kappa}$  
to define the co-product in (\ref{newco}). Thinking of 
the $\bk_1,\bk_2$ as momenta
we should think of (\ref{BCH}) as a ``curved'' momentum addition,
which reduces to the usual linear 
 addition in the limit $\kappa\rightarrow 0$.

The universal $R$-element (\ref{bigrmatrix})
  of  $D(SU(2))\otimes D(SU(2))$  can be mapped  into $A_\kappa\otimes
A_\kappa$
 via the pull-back  $\hbox{EXP}^*_\kappa\otimes\hbox{EXP}^*_\kappa$.
We denote the pull-back by $R_{\kappa}$, i.e.
\bea
\label{rkappa}
R_{\kappa}=\bigl(\hbox{EXP}^*_\kappa\otimes\hbox{EXP}^*_\kappa\bigr)( R ).
\eea    
 Thus we obtain a family
 $(A_{\kappa},\bullet,\Delta_\kappa,1,\epsilon,S,R_{\kappa}$)
of quasi-triangular Hopf algebras which are all isomorphic 
to $D(SU(2))$ for $\kappa >0$.
The last remark in the previous paragraph 
 shows that the co-multiplication in  $A_\kappa$ tends
to the co-commutative co-product of   $A_0$ (\ref{isocomult})
in the limit $\kappa\rightarrow 0$.

The relation between
the universal $R$-element $R_\kappa$ 
 and the $r$-matrix (\ref{rmatrix})  of $iso(3)$  can now be
established as follows.
Since  $R_\kappa \in A_\kappa\otimes A_\kappa \subset A_0 \otimes A_0$
we can let it act in an  irrep $\Pi_{M_1,s_1}\otimes \Pi_{M_2,s_2}$
 of $ A_0\otimes A_0$ and 
consider the limit
$\kappa\rightarrow 0$. To keep the following formulae simple
 we write $\Pi_1 \otimes \Pi_2$ instead of 
$\Pi_{M_1s_1}\otimes \Pi_{M_2s_2}$.
Then
\bea
R_\kappa (\bk_1,y_1,\bk_2,y_2) =
\delta_e(\exp_\kappa(\bk_1)y_2^{-1})\delta_e(y_1).
\eea
Acting on $\Phi\in H_{s_1}\otimes H_{s_2}$ via 
$\Pi_1\otimes \Pi_2$ gives 
\bea
\bigl(  \Pi_1\otimes \Pi_2
(R_\kappa)
\Phi\bigr)(x_1,x_2) = \Phi(x_1,\exp_{\kappa}\bigl(\bk(x_1,M_1)\bigr)\,x_2),
\eea
where we have used the notation (\ref{bpdef}).
Then we find for small  $\kappa$ 
\bea
\bigl(\Pi_1\otimes \Pi_2
(R_\kappa)\,\Phi\bigr) (x_1,x_2) = 
\bigl( 1\otimes 1 +  \kappa\,
 k^a(x_1,M_1) {\mathcal J}^{(2)}_a \bigr) \,\Phi(x_1,x_2)\,+\,{\mathcal O}
(\kappa^2).
\eea
The term of order  $\kappa$  is  the action of $r$ in the representation
$\Pi_1\otimes \Pi_2$ (\ref{raction}).
Since the above statement is true in any representation, we conclude
that  in the limit  $\kappa\rightarrow 0$
\bea
\label{result}
R_\kappa = 1\otimes 1 + i \kappa\,  r + \,{\mathcal O}(\kappa^2).
\eea

\section{Discussion and outlook}

In the previous section we have implemented the key step of 
the combinatorial quantisation procedure for the case of  Euclidean
 three dimensional gravity. The quantum $R$-matrix $R_{\kappa}$ 
can now be used to define the graph algebra as introduced and  discussed  in 
\cite{AGSI} and \cite{AGSII}.
The relations defining this algebra are numerous
and complicated and will not be repeated here.  From the graph 
algebra one obtains another algebra, introduced in \cite{AGSII} 
and  called  the moduli algebra ${\mathcal M}$ in  \cite{AS},  
by  taking a quotient
and imposing various conditions.
The moduli algebra is the   quantisation
of the commutative algebra of functions on the classical moduli
space $M(G,g,\calC_1,...,\calC_m)$, which 
in our case is interpreted as  the 
 gravitational phase space. The algbra 
${\mathcal M}$ may thus be thought of as  
the algebra of observables in quantised three-dimensional gravity.
It would clearly be very interesting to study this algebra
in more detail, and to interpret it physically. This should 
probably be done in the case of Lorentzian gravity, where the physical
interpretation is clearer. One important aspect of that 
discussion will be a clear understanding of the classical
limit. Here we only note that this limit is related to
the limit $\kappa \rightarrow 0$ of our deformation parameter, 
as can be seen by a simple dimensional analysis.
 In units where the speed of light $c$ is 1, the 
gravitational coupling constant $G_3$ in three dimensions (Newton's
constant) has the dimension  of inverse mass. The elements $\bk$ 
of the space  $(\RR^3)^*$ defined in the previous section have the dimension
of inverse length. Since the argument of the 
 exponential function (\ref{expfun})
must be dimensionless, it follows that the deformation parameter
$\kappa$ has the dimension of  length. Since $G_3$ and 
Planck's constant $\hbar$ are the only available physical constants
in pure three dimensional quantum gravity (having set $c=1$), 
it follows that $\kappa$ is proportional to $\hbar G_3$. The
classical limit $\hbar\rightarrow 0$ therefore corresponds, for 
fixed $G_3$, to the limit  $\kappa\rightarrow 0$.

The final step in the quantisation program, carried 
out for compact gauge groups in \cite{AS},  is the construction of 
a representation of the moduli  algebra ${\mathcal M}$ 
on a suitable Hilbert space. 
The implementation of this step for three dimensional quantum
gravity would again be of great interest. For Riemann surfaces
with genus zero one almost certainly reproduces the results  of 
\cite{BM}, where  the Hilbert space 
of $m$ distinguishable particles on a genus zero surface 
is constructed from irreps of $D(SU(2))$.
 As pointed out in \cite{BM} the $R$-element
of $D(SU(2))$ provides a representation of the braid group action
on the multi-particle Hilbert space. As further explained in
\cite{Muller}, the analogous   braid group action 
in the Lorentzian case 
 allows one 
to compute  the gravitational quantum  scattering of particles,
at least in principle. In practice there are a number of conceptual
and practical problems, addressed  in  \cite{BMS}.
Within the combinatorial quantisation programme
it is not difficult also to  consider  particles on 
Riemann surfaces of higher genus.
As explained in \cite{AS} 
 one is guaranteed to obtain a Hilbert space which carries a 
representation  of the mapping class group of the surface.
In gravity, the requirement of 
diffeomorphism invariance implies 
that only states which are invariant
under the action of the mapping class group are phyiscal.
The scheme outlined here thus leads to a natural
 implementation of diffeomorphism invariance in quantum gravity,
something which has proved difficult in other approaches.

Mathematically, a number of issues 
 need to be resolved.
The classical moduli spaces $M(G,g,\calC_1,...,\calC_m)$ 
for the non-compact gauge groups in table 1 should be investigated
rigorously, and both their symplectic and their   
singularity structure should be clarified. Here the more 
geometric approach to the symplectic structure developed in
\cite{GHJW} and \cite{AMM} may be useful.
Looking ahead, a number of generalisations of our sects. 4 and 5
should be studied.
The deformation of the group algebra  of $ISO(3)^\sim$
to the quantum double $D(SU(2))$  described in sect. 5 of this paper
can be generalised to 
 semi-direct product groups of the form $\gothh\rtimes H$, where
$H$ is an arbitrary locally compact unimodular Lie group 
and  $\gothh$ the Lie algebra of $H$,
thought of as an abelian group under addition, with 
$H$ acting  via conjugation. Applying the procedure described
in this paper, 
the group algebra of $\gothh\rtimes H$   can be deformed to 
$D(H)$. More ambitiously,  it would be interesting to include
the cosmological constant in our discussion.  
In \cite{BR} the quantum group 
$SL_q(2,\CC)_\RR$ was argued to be relevant in the 
combinatorial quantistation of $SL(2,\CC)$ Chern-Simons theory.
This would suggest an interesting connection between  the quantum group 
$SL_q(2,\CC)_\RR$ and the quantum
 doubles $D(SU(2))$ and $D(SO(2,1))$.

\vspace{2cm}

\vbox{
\noindent{\bf Acknowledgments}

\noindent
The question addressed
in this article arose in discussions I had
 with Sander Bais and Nathalie Muller while I was a post doc
in the Institute for Theoretical Physics of the University of 
Amsterdam. Since then I have had further illuminating discussions
with  Anton Alekseev, Sander Bais, David Calderbank,   Tom Koornwinder, Volker Schomerus,
  Allan Sinclair and Joost Slingerland. I thank the organisers of 
the workshop ``Deformation quantisation of singular reduced spaces''
in Oberwolfach for the opportunity to present an 
 early version of this work, and 
 acknowledge financial support through an Advanced Research
Fellowship of the Engineering and Physical Sciences Research Council.

}


\begin{thebibliography}{99}

\itemsep=\smallskipamount
\bibitem{AGSI}
A.~Y.~Alekseev, H.~Grosse and V.~Schomerus, {\em Combinatorial
quantization of the Hamiltonian Chern-Simons Theory,}
Commun.~Math.~Phys., {\bf 172} (1995), 317--358.
\bibitem{AGSII}
A.~Yu.~Alekseev, H.~Grosse and V.~Schomerus, {\em Combinatorial
quantization of the Hamiltonian Chern-Simons Theory II,}
Commun.~Math.~Phys., {\bf 174} (1995), 561--604.
\bibitem{AS} A.~Yu.~Alekseev and V.~Schomerus, {\em Representation 
theory of Chern-Simons observables,} Duke Math.~Journal, {\bf 85}
(1996), 447--510.
\bibitem{FR} V.~V.~Fock and A.~A.~Rosly, {\em Poisson structures on
moduli of flat connections on Riemann surfaces and $r$-matrices,} 
 ITEP preprint {\bf  72-92}  (1992); see also {\tt math.QA/9802054}.
\bibitem{BM} F.~A.~Bais and N.~M.~Muller, {\em  Topological field
 theory and the quantum double of $SU(2)$,}
Nucl.~Phys., {\bf  B530} (1998), 349--400.
\bibitem{BMS} F.~A.~Bais, N.~M.~Muller and B.~J.~Schroers, {\em
Quantum double symmetry and topological interactions in
(2+1)-dimensional  quantum gravity}, in preparation.
\bibitem{BR} E.~Buffenoir and Ph.~Roche, {\em Harmonic analysis on the quantum
Lorentz group,} Commun.~Math.~Phys., {\bf 207} (1999), 499--555.
\bibitem{AT} A.~Achucarro and P.~Townsend, {\em  A Chern--Simons
 action for three-dimensional anti-de Sitter supergravity
 theories,}  Phys.~Lett., {\bf B180} (1986),  85--100.
\bibitem{Witten1} E.~Witten, {\em  2+1 dimensional gravity as an exactly
 soluble system,} Nucl.~Phys.,  {\bf B311} (1988), 46--78.
\bibitem{Carlipbook} S.~Carlip, {\em Quantum gravity in (2+1) dimensions,}
Cambridge University Press, Cambridge, {\bf 1998}.
\bibitem{Carlipart} S.~Carlip, {\em Lectures on 2+1 dimensional gravity,}
UCD-95-{\bf 6} (1995), {\tt gr-qc/9503024}.
\bibitem{Sharpe} R.~W.~Sharpe, {\em Differential Geometry},
Springer Verlag, New York, {\bf 1996}.
\bibitem{Matschull}
H.-J.~Matschull,
{\em On the relation between (2+1) Einstein gravity and Chern-Simons Theory,}
Class.~Quant.~Grav., {\bf 16} (1999), 2599--2609. 
\bibitem{MR}
J.~E.~Marsden and T.~S.~Ratiu, {\em Introduction to mechanics and symmetry,}
 Springer Verlag, New York, {\bf 1994}.
\bibitem{Witten2} E.~Witten, {\em Quantization of Chern-Simons gauge
theory with complex gauge group,} Commun.~Math.~Phys., {\bf 137}
(1991), 29--66.
\bibitem{AB} M.~Atiyah and R.~Bott, {\em The Yang-Mills equations over
Riemann surfaces,} Phil.~Trans.~Roy.~Soc. London Ser. A, {\bf 308}
(1983), 523--615.
\bibitem{Atiyah} M.~Atiyah, {\em The geometry and physics of knots,}
Cambridge University Press, Cambridge, {\bf 1990}.
\bibitem{CP} V.~Chari and A.~Pressley {\em Quantum Groups, }
Cambridge University Press,  Cambridge, {\bf 1994}.
\bibitem{Stachura}
 P.~Stachura, {\em Poisson-Lie structures on 
Poincar\'e and Euclidean groups in three dimensions,} J.
J.~Phys., {\bf A 31} (1998), 4555--4564.
\bibitem{CG}
E.~Celeghini and R.~Giachetti, {\em The three dimensional euclidean
quantum group $E(3)_q$ and its R-matrix,} J.~Math.~Phys., {\bf 32}
(1991), 1159--1165. 
\bibitem{DVVV}
R.~Dijkgraaf, C.~Vafa, E.~Verlinde, and H.~Verlinde,
{\em The operator algebra of orbifold models,}
Commun.~Math.~Phys., {\bf 123} (1989), 485--526.
\bibitem{BDP}
F.~A.~Bais, P.~van~Driel, and M.~de~Wild~Propitius,
{\em Quantum symmetries in discrete gauge theories,}
Phys.~Lett., {\bf B 280} (1992), 63--70.
\bibitem{DPR} R.~Dijkgraaf, V.~Pasquier, and  P.~Roche, 
 {\em Quasi Hopf algebras, group cohomology and orbifold models,} 
Nucl.~Phys. (Proc. Suppl.), {\bf 18B} (1990), 60--72.
\bibitem{Bonneau}
P.~Bonneau, {\em Topological quantum double,}  Rev.~Math.~Phys., {\bf 6}
(1994), 305--318.
\bibitem{Mueger}
M.~M\"uger, {\em Quantum Double actions on operator algebras and 
orbifold quantum field theories,} Commun.~Math.~Phys., {\bf 191} (1998),
137--181.
\bibitem{KM}
T.~H.~Koornwinder and N.~M.~Muller.
{\em The quantum double of a (locally) compact group,}
J.~Lie Theory {\bf  7} (1997), 33--52;{\bf 8} (1998), 187 (erratum).
\bibitem{KBM} T.~H.~Koornwinder, F.~A.~Bais and N.~M.~Muller,
{\em Tensor Product Representations of the Quantum Double of a Compact
Group,} Commun.~Math.~Phys., {\bf 198} (1998), 157--186.
\bibitem{Pedersen}
G.~K.~Pedersen, $C^*${\em algebras and their automorphism groups,}
Academic Press, London, {\bf 1979}. 
\bibitem{GHJW}
K.~Guruprasad, J.~Huebschmann, L.~Jeffrey, A.~Weinstein, {\em 
Group systems, groupoids, and moduli spaces of parabolic bundles,} 
Duke Math.~J., {\bf 89 } (1997), 377--412. 
\bibitem{AMM}
A.~Yu.~Alekseev, A.~Z.~Malkin and E.~Meinrenken, {\em Lie group
valued moment maps,} J.~Differential Geom., {\bf 48} (1998), 445--495.
\bibitem{Muller}
N.~Muller, {\em Topological interactions and quantum double symmetries,}
Ph.D. dissertation, University of Amsterdam, {\bf 1998}.
\end{thebibliography}
\end{document}